\documentclass[12pt,a4paper]{amsart}
\usepackage{amsmath, amsfonts, xifthen, latexsym, amssymb, amsthm, amscd}
\usepackage[utf8]{inputenc}
\usepackage{graphicx}
\usepackage{url}
\usepackage{setspace}
\usepackage{epigraph}
\theoremstyle{plain}
\newtheorem{theorem}{Theorem} \newtheorem{lemma}{Lemma}[section]
\newtheorem{propo}{Proposition}[section]

 \newtheorem{defin}{Definition}[section]

 \newcommand{\eps}{\epsilon} 
 
\newcommand{\R}{\mathbb{R}}

\newcommand{\cH}{\mathcal{H}}

\newcommand{\cG}{\mathcal{G}}
\newcommand{\uc}{\underline{c}}

\newcommand{\vi} {\vskip 0.1in \noindent}
\title{Uniform Local Amenability implies Property A}
\author{G\'abor Elek}
\address{Department of Mathematics And Statistics, Fylde College, Lancaster University, Lancaster, LA1 4YF, United Kingdom
\vskip 0.05in and
\vskip 0.1in Alfred Renyi Institute, Budapest, Hungary}
\email{g.elek@lancaster.ac.uk}  
\thanks{The author was partially supported
by the ERC Starting Grant "Limits of Structures in Algebra and Combinatorics"  }
\begin{document}
\begin{abstract}
In this short note we answer a query 
of Brodzki, Niblo, \v{S}pakula, Willett and Wright \cite{ULA} by showing that
all bounded degree uniformly locally amenable graphs have Property A. For the second result of the note
recall that Kaiser \cite{Kaiser} proved that if $\Gamma$ is a finitely generated group and $\{H_i\}^\infty_{i=1}$ is
a Farber sequence of finite index subgroups, then the associated Schreier graph sequence is of Property A if and only if the group is amenable.
We show however, that there exist a non-amenable group and a nested sequence of finite index subgroups $\{H_i\}^\infty_{i=1}$ such that
$\cap H_i=\{e_\Gamma\}$, and the associated Schreier graph sequence is of Property A. 
\end{abstract}
\maketitle
\section{Introduction}
First, let us introduce the main notions of our paper: Local Hyperfiniteness,
Uniform Local Amenability, Property A and Local Strong Hyperfiniteness,
\subsection{Local Hyperfiniteness}
First of all, fix an integer $d>0$.
In the course of the note,  all graphs, finite or infinite, are considered to be of vertex degree bound $d$. We denote the class of
all finite graphs of vertex degree bound $d$ by $\cH_d$. 
A class of finite  graphs $\cG\subset\cH_d$ is called {\bf hyperfinite} (see \cite{Elekcost}) if for any $\epsilon>0$, there exists some $K>0$
such that from any graph $G\in\cG$ one can remove at most $\eps |V(G)|$ vertices (with all the incident edges) in such a way
that all the components of the resulting graph $G'$ have at most $K$ vertices. Note that the class of planar graphs 
or graphs of fixed subexponential growth are hyperfinite.
\begin{defin} A not necessarily connected
infinite graph $S$ is called {\bf locally hyperfinite} if the class $\cG_S$ of all finite subgraphs of $S$ is
a hyperfinite class. \end{defin}
\subsection{Uniform Local Amenability}
The notion of Uniform Local Amenability was introduced by Brodzki, Niblo, \v{S}pakula, Willett and Wright. Since we are interested only in bounded degree
graphs, the definition of Uniform Local Amenability fits into our scheme of notions rather nicely.
A class of finite graphs $\cG$ is called {\bf uniformly locally amenable} if for any $\eps>0$ there exists some $K>0$ such that
for any $G\in\cG$ one has a subset $E\subset V(G)$, $|E|\leq K$ such that $|\partial (E)|\leq \eps|E|$. Here, $\partial(E)$ denotes the set of vertices
in $E$ which are adjacent to some vertices that are not in $E$. 
\begin{defin}
 An infinite graph $S$ is called {\bf uniformly locally amenable} if the class $\cG_S$ is uniformly locally amenable.
\end{defin} 
\noindent
It is interesting to note that Lemma 2.4 of \cite{ElekL2} the uniform local amenability of graph classes of subexponential growth
was observed (in this paper hyperfiniteness was called ``antiexpander property'', the current terminology appeared first in \cite{Elekcost}.
\subsection{Property A}
Property A was introduced by Guoliang Yu \cite{Yu}. Instead of the original definition, we will use the notion
of weighted hyperfiniteness \cite{Elektimar}, since Sako \cite{Sako} proved the equivalence of the two notions.
Let $G$ be a finite graph. We say that $G$ is weighted $(\eps,K)$-hyperfinite, if for any nonnegative function
$w:V(G)\to \R$ there exists a subset $Y\subset V(G)$ such that
\begin{itemize}
\item $\sum_{y\in Y} w(y)\leq \eps \sum_{x\in V(G)} w(x)$,
\item if we remove $Y$ from $V(G)$ with all the incident edges, then the size of the components of the resulting graph
are at most $K$.
\end{itemize}
\begin{defin}
An infinite graph $S$ is of Property A (that is, Locally Weighted Hyperfinite) if the class $\cG_S$ is weighted hyperfinite.
\end{defin}
\noindent
One should note that the Cayley graph of a finitely generated group $\Gamma$ is of Property A if and only if $\Gamma$ is exact. 
Hyperbolic groups, amenable groups are exact. Directed unions, subgroups and extensions of exact groups are also exact.
However, in 2003 Gromov \cite{Gromov} constructed a non-exact group, so the Cayley graph of his group is not of Property A.
\subsection{Local Strong Hyperfiniteness}
The notion of Strong Hyperfiniteness was introduced in \cite{Elek2} (under the name of ``uniform hyperfiniteness''). 
Let $G$ be a finite graph and $Y\subset V(G)$
is a subset of vertices. We say that $Y$ is a $K$-separator if we remove $Y$ from $V(G)$ with all the incident edges, then the size of the components of the resulting graph
are at most $K$.
We say that a graph $G$ is strongly $(\eps,K)$-hyperfinite if there is a probability measure $\mu$ on the set of $K$-separators
such that for all $x\in V(G)$ the measure of $K$-separators containing $x$ is less or equal than $\eps$.
\noindent
So, a class of finite graphs $\cG$ is called {\bf strongly  hyperfinite} if for any $\eps>0$ there exists a $K>0$ such that
all $G\in\cG$ is strongly $(\eps,K)$-hyperfinite.
\begin{defin} An infinite graph $S$ is called {\bf locally strongly hyperfinite} if the class $\cG_S$ is strongly
hyperfinite. \end{defin}
Now, the plot thickens. The notion of {\bf fractional cc-fragility} for graph classes was defined by Romero, Wrochna and \v{Z}ivn\'y \cite{RWZ}, which
is clearly equivalent to strong hyperfiniteness. 
One should note that the idea of fragility goes back to Dvorak \cite{Dvo}.
The results of Romero, Wrochna and \v{Z}ivn\'y will be absolutely crucial for our note.
\subsection{The main theorem}
It was proved in \cite{ULA} that Property A implies Uniform Local Amenability. The authors asked whether the converse is true. 
Our main result is the following theorem.
\begin{theorem}\label{main}
If $G$ is a uniformly locally amenable infinite graph of bounded vertex degrees, then $G$ is of Property A.
\end{theorem}
\noindent
Let $\Gamma$ be a finitely generated group with a finite, symmetric generating system $\Sigma$. For a finite index subgroup $H$, we can consider the left Schreier graph
$\mbox{Sch}(\Gamma/H,\Sigma)$.
Kaiser \cite{Kaiser} proved that if $\Gamma$ is a finitely generated group and $\Gamma\supset H_1 \supset H_2 \supset\dots,\,\cap_{i=1}^\infty H_i=\{e_\Gamma\}$ is
a Farber sequence of finite index subgroups, then the associated Schreier graph sequence is of Property A if and only if $\Gamma$ is amenable.
Using our results we will see that the Farber Property is, in fact, crucial, since we have the following proposition.
\begin{propo} \label{Kai}
There exist a finitely generated non-amenable group and finite index subgroups $\Gamma\supset H_1 \supset H_2 \supset\dots,\,\cap_{i=1}^\infty H_i=\{e_\Gamma\}$ 
such that the associated Schreier graph sequence is of Property A.
\end{propo} 
\vskip 0.1in
\noindent
{\bf Acknowledgement:} The author would like to thank Ana Khukhro for calling his attention to the paper of Kaiser. Also, we are indebted to Endre Cs\'oka for
providing us with the very short proof for Proposition \ref{weight}.

\section{Uniform Local Amenability implies Local Hyperfiniteness} \label{UiL}
\begin{propo} \label{elso}
Every uniformly locally amenable infinite graph $S$ is locally hyperfinite.
\end{propo}
\proof It is enough to prove that any monotone (that is closed under taking subgraphs) uniformly locally amenable graph class is hyperfinite. We just repeat the argument 
of Lemma 2.4 in \cite{ElekL2}.  Let $\cG$ be a monotone uniformly locally amenable class of graphs and $\eps>0$ and $K>0$ be numbers as in 
the definition of uniform local amenability. Let $G\in \cG$. Find $E_1\subset V(G)$ such that $|\partial_G(E_1)|\leq \eps |E_1|$ and $|E_1|\leq K$. Remove $E_1$ from $G$ to obtain
the graph $G_2$. By monotonicity, we have $E_2\subset V(G)$ such  $|\partial_{G_2}(E_2)|\leq \eps |E_2|$ and $|E_2|\leq K$. Now we remove $E_2$. Inductively, we construct subsets $E_3,E_4,\dots$ and
subgraphs
$G_3,G_4,\dots$ 
in such a way that $|\partial_{G_i}(E_i)|\leq \eps |E_i|.$ Hence $\cup \partial_{G_i}(E_i)$ will be an $(\eps,K)$-separator for $G$. \qed.

\section{Local Hyperfiniteness implies Local Strong Hyperfiniteness}
By Lemma 8.1 of \cite{RWZ}, every hyperfinite monotone graph class is strongly hyperfinite. Since for an infinite graph $S$, the graph class $\cG_S$ is monotone 
we immediately have the following proposition.
\begin{propo} \label{masodik}
Every locally hyperfinite infinite graph $S$ is locally strongly hyperfinite.
\end{propo}

\section{Local Strong Hyperfiniteness implies Property A}
\begin{lemma} A finite graph $G$ is $(\eps,K)$-strongly hyperfinite if and only if it is $(\eps,K)$-weighted hyperfinite as well. \label{weight}
\end{lemma}
\proof
Let $G$ be a finite graph and $V(G)=n$. Let $k$ be the number of $K$-separators in $G$. For such a $K$-separator $Y$ let $\{\underline{c}_Y:V(G)\to \{0,1\}\}\in \R^n$ its characteristic vector.
We define the hull $\cH$ of the $K$-separators as the convex set of vectors $\underline{z}\in \R^n$ which can be written in the form
$$\underline{z}=\sum_{i=1}^k x_i \underline{c}_{Y_i}+\underline{y}\,,$$
\noindent
where all $x_i\geq 0$, $\sum^k_{i=1} x_i=1$ and $\underline{y}$ is a non-negative vector.
Let $\underline{v}=\{\eps,\eps,\dots,\eps\}$.
Note that $G$ is $(\eps,K)$-strongly hyperfinite 
if and only if
there exists non-negative real numbers $\{x_i\}^k_{i=1},
\sum^k_{i=1} x_i=1$ such that
$$\sum_{i=1}^k x_i \uc_{Y_i}\leq \underline{v}\,.$$
\noindent
Also, $G$ is $(\eps,K)$-weighted hyperfinite
if and only if for any non-negative vector $\underline{w}$ of total weight $1$,  there
exists a $K$-separator $Y$ such that
$$\langle \underline{w}, \uc_Y \rangle \leq \eps\,.$$
\noindent 
We have two cases.
\vi
{\bf Case 1.} $\underline{v}\in \cH$, that is, 
$G$ is $(\eps,K)$-strongly hyperfinite.
Let $\underline{w}\in \R^n$ be a non-negative vector of total weight $1$.
Suppose that $\langle \underline{w}, c_Y \rangle >\eps$
holds for all the $K$-separators $Y$.
Then, $\langle \underline{w}, \underline{h} \rangle >\eps$
holds for all $h\in \cH$. Hence, we have
$$\eps=\langle \underline{w}, \underline{v} \rangle >\eps\,,$$
\noindent
leading to a contradiction.
Hence, there exists a $K$-separator $Y$ such that
$\langle \underline{w}, c_Y \rangle \leq \eps\,.$
So, $G$ is $(\eps,K)$-weighted hyperfinite.

\noindent
{\bf Case 2.} $\underline{v}\notin \cH$. That is,
$G$ is not $(\eps,K)$-strongly hyperfinite. Since $\cH$ is a closed convex set, there must exist a hyperplane $H\in \R^n$ containing
$\underline{v}$ such that $\cH$ is entirely on one side of the hyperplane $H$.
That is, there exists a vector $\underline{w}\in \R^n$ such that for any $\underline{h}\in \cH$ we have
$$\langle \underline{w},\underline{v} \rangle < \langle \underline{w},\underline{h} \rangle\,.$$
\noindent
Then, for any $1\leq i \leq n$, $w_i\geq 0$. Indeed,
if there exists $1\leq i \leq n$ such that $w_i<0$, then for any vector $c_Y$ there exists
a non-negative vector $\underline{y}\in \R^n$ such that
$$\langle \underline{w},\uc_Y+\underline{y} \rangle < \langle \underline{w},
\underline{v} \rangle \,.$$
\noindent
We can also assume that the total weight
of $\underline{w}$ is $1$.
So, 
$\eps< \langle \underline{w}, \underline{c}_Y \rangle $
holds for each $K$-separator $Y$, that is, $G$ is not $(\eps, K)$-weighted hyperfinite. \qed
\vi
Hence, we have the following proposition.
\begin{propo} \label{harmadik}
Every locally strongly hyperfinite infinite graph $S$ is of Property A.
\end{propo}
\noindent
Now, our Theorem follows from Proposition \ref{elso}, Proposition \ref{masodik} and Proposition \ref{harmadik}.

\section{A remark on a theorem of Kaiser}
In \cite{Abert}, Mikl\'os Ab\'ert and the author considered the non-amenable group $\Gamma=C_2*C_2*C_2*C_2$, the free product of four cyclic graphs with its
canonical generating system of four elements and constructed
a sequence of finite index subgroups $\{H_i\}^\infty_{i=1}$ such that
\begin{itemize}
\item $\Gamma\supset H_1 \supset H_2 \supset \dots$
\item $\cap_{i=1} H_i=\{e_\Gamma\}.$
\item Each Schreier graph $\mbox{Sch}(\Gamma/H_i,\Sigma)$ is a planar graph.
\end{itemize}
\noindent
Since planar graphs of bounded vertex degrees form a monotone hyperfinite class, by our main theorem, the system of the Schreier graphs above is of Property A. 
Hence, Proposition \ref{Kai} immediately follows. \qed

\end{document}